\documentstyle[12pt]{article}

  \newcommand{\const}{\rm const}
  \newcommand{\Var}{\rm Var}
  \newcommand{\supp}{\rm supp}

\begin{document}

   \begin{center}

  \   {\bf   Bilateral tail estimate for distribution of self normalizes }\\

\vspace{4mm}

   \  {\bf   sums of independent centered random variables }\\

\vspace{4mm}

 \ {\bf  under natural norming }\\

\vspace{4mm}

  \  {\bf Ostrovsky E., Sirota L.}\\

\vspace{4mm}

 Israel,  Bar-Ilan University, department of Mathematic and Statistics, 59200, \\

\vspace{4mm}

E-mails: eugostrovsky@list.ru, \hspace{5mm} sirota3@bezeqint.net \\

\vspace{5mm}

  {\bf Abstract} \\

\vspace{4mm}

 \end{center}

\ We derive in this article the exact non-asymptotical exponential and power estimates for self-normalized sums
of centered independent random variables (r.v.) under natural norming. \par

 \ We will use also the theory of the so-called Grand Lebesgue Spaces (GLS) of random variables.\\

\vspace{5mm}

 \par

\vspace{4mm}


\vspace{5mm}

\section{ Definitions.  Notations. Previous results.  Statement of problem.}
 \vspace{3mm}

 \ Let  $ \ \{\xi(i) \}, \ i = 1,2,\ldots,n; \ \xi := \xi(1)  $ be a sequence of centered: $ {\bf E} \xi(i) = 0 $
independent identically distributed (i., i.d.) random variables (r.v.) defined on
certain probability space, having a finite non-zero variance $  \sigma^2 := {\bf E} \xi^2(i) \in (0, \infty). $
Introduce the following self-normalized sequence of sums under natural norming

$$
T(n) = \sqrt{n} \cdot \frac{\sum_i \xi(i)}{\sum_i \xi^2(i)}, \eqno(1.1)
$$
here and in what follow

$$
\sum = \sum_i = \sum_{i=1}^n,
$$
and define the correspondent tail probabilities

$$
Q_n = Q_n(B) := {\bf P}(T(n) > B), \ B = \const > 0; \ Q = Q(B) := \sup_n Q_n(B). \eqno(1.2)
$$

\vspace{4mm}

  \ {\bf Our purpose in this preprint is obtaining non-asymptotical exponential and power bounds for introduced in (1.2) tail probabilities. } \par

\vspace{4mm}

 \ This problem with another self norming sequence was considered in many works,
see e.g. [3], [7], [8], [10], [16], [17], [18]-[19], [21]. Note that in these works
was considered as a rule only asymptotical approach, i.e. when $ \ n \to \infty; \ $  for
instance, was investigated the classical Central Limit Theorem (CLT), Law of
Iterated Logarithm (LIL) and Large Deviations (LD) for these variables.
Several interest applications in the non-parametrical statistics are described in [3], [7], [17], [18]-[19] etc. \par

\vspace{4mm}

 \ We must introduce now some needed notions and notations. $ \ \eta(i) = \eta(i; n, B)  :=  $

$$
 \sqrt{n} \ \xi(i) + B (\sigma^2 - \xi^2(i)),  \ \eta = \eta(n; B)  = \eta(1) = \sqrt{n} \ \xi + B \ (\sigma^2 - \xi^2);
$$

$$
\gamma(i) = \gamma(i; \ n, \ B) := \eta(i)/\sqrt{n} = \xi(i)  + B  (\sigma^2 - \xi^2(i)) /\sqrt{n}, \ \gamma = \gamma(1);
$$

$$
S = S(n) = \sum \xi(i), \  V = V(n) = \sum (\xi^2(i) - \sigma^2), \eqno(1.3)
$$
so that

$$
{\bf E} S(n) = {\bf E} V(n) = 0
$$
and

$$
T(n) = \frac{ \sqrt{n} \cdot S(n)}{V(n) + n \ \sigma^2}. \eqno(1.4)
$$

 \ Further, we introduce the functions {\it of two variables}

$$
\phi(\lambda_1, \lambda_2) \stackrel{def}{=} \ln {\bf E} \exp ( \lambda_1 \xi + \lambda_2 (\sigma^2 - \xi^2)     ),   \eqno(1.5)
$$
so that

$$
{\bf E} \exp \left( \mu_1 \ \sqrt{n} \xi  + \mu_2 \ B \ (\sigma^2 - \xi^2)  \right) = \exp \phi( \mu_1 \sqrt{n}, \mu_2 \ B)  ).  \eqno(1.6)
$$

 \ Denote as ordinary for any r.v.  $ \ \zeta \ $ its classical Lebesgue-Riesz $  \ L(p) \ $ norm

$$
|\zeta|_p := [ {\bf E} |\zeta|^p  ]^{1/p}, \ p \ge 1,
$$
and introduce the variables

$$
m(p) = m(p; \ B, \sigma, \ n) := | n^{-1/2} \ \xi  +  B \ (\sigma^2 - \xi^2)  |_p = n^{-1/2} |\eta|_p, \eqno(1.7)
$$
if of course  $  m(p) \ $  is finite for certain value $  \ p, \ p \ge 1; \ $

$$
w = w(\sigma) = w(\sigma; \xi) := {\bf E} (\sigma^2 - \xi^2)^2; \ z = z(\sigma) = z(\sigma; \xi) := {\bf E} (\sigma \xi - \xi^3), \eqno(1.8)
$$
so that

$$
D^2(\sigma, n, B; \xi) := \Var(\eta) = n \ \sigma^2 + 2 \ B \ \sqrt{n} \ z + B^2 \ w \eqno(1.9)
$$
and

$$
{\bf E} e^{\theta  \ \eta(i;n)} = e^{ \nu(\theta) } = e^{  \nu(\theta; n, B) }, \eqno(1.10)
$$

where

$$
 \nu(\theta) = \nu(\theta; n, B)   =   \phi(\theta \ \sqrt{n}, \ B \ \theta). \eqno(1.11)
$$

\vspace{4mm}

 \section{ Grand Lebesgue Spaces (GLS).}

\vspace{4mm}

  \ Let  $  Z = (Z, M, \mu) $ be probability space with non-trivial measure  $ \mu. $ Let
also $  \psi = \psi(p), \ p \in [1,b), \ b = \const \in (1, \infty] $ be certain bounded
from below: $ \inf \psi(p) > 0 $ continuous inside the semi-open interval  $ p \in [1,b) $
numerical function. We can and will suppose  $  \ b = \sup \{ p, \ \psi(p) < \infty  \}, $
 so that   $ \supp \psi = [1, b) $ or  $ \supp \psi = [1, b]. $ The set of all such a functions will be
denoted by  $ \ \Psi(b) = \{ \psi(\cdot)  \}; \ \Psi := \Psi(\infty).  $ \par

\vspace{4mm}

 \ By definition, the (Banach) Grand Lebesgue Space (GLS) space   $  \ G \psi  = G\psi(b) $
consists on all the real (or complex) numerical valued measurable functions
(random variables, r.v.)   $   \ \zeta \ $  defined on our probability space and having a finite norm

$$
||\zeta|| = ||\zeta||G\psi \stackrel{def}{=} \sup_{p \in [1,b)} \left[ \frac{|\zeta|_p}{\psi(p)} \right]. \eqno(2.1)
$$

Here $ \  |\zeta|_p = |\zeta|L_p(Z).  \  $ \par

 \ These spaces are Banach functional space, are complete, and rearrangement
invariant in the classical sense, see [1], chapters 1, 2; and were investigated in
particular in many works, see e.g. [2], [5], [6], [9], [11], [12]-[15] etc. We refer
here some used in the sequel facts about these spaces and supplement more. \par

 \ It is known that if  $  \zeta \ne 0, $ then

$$
{\bf P} (|\zeta| > y) \le  \exp \left( -v_{\psi}^*(\ln(u/||\zeta||)   \right). \eqno(2.2)
$$

 \ Here and in the sequel the operator  $ f \to f^*  $   will denote the Young-Fenchel transform

$$
f^*(u) \stackrel{def}{=} \sup_{x \in Dom(f)} (x u - f(x)).
$$

 \ Conversely, the last inequality may be reversed in the following version: if

$$
{\bf P} (\zeta > u) \le  \exp \left( - v_{\psi}^* (\ln (u/K) \right), \ u \ge e \ K,
$$
and if the auxiliary function

$$
v(p) = v_{\psi}(p)  \stackrel{def}{= } p \ln \psi(p), \ p \in [1,b)
$$
is positive, continuous, convex and such that

$$
\lim_{p \to \infty} \psi(p) = \infty,
$$
then $  \ \zeta \in G(\psi) \  $ and besides $  \ ||\zeta|| \le C(\psi) \cdot K.  $\par

\vspace{4mm}

 \ Let us consider the so-called exponential Orlicz space $  L(M) $ builded over
source probability space with correspondent Young-Orlicz function

$$
M(y) = \exp \left(  v_{\psi}^*(\ln y), \right) \ y \ge e;  \ M(y) = C y^2, \ |y| < e.
$$

 \ The Orlicz $  ||\cdot||  $ and $ G\psi $ norms are quite equivalent:

$$
||\zeta||G\psi  \le C_1 ||\zeta||L(M) \le C_2 ||\zeta||G\psi, \ 0 < C_1 < C_2 < \infty. \eqno(2.3)
$$

\vspace{4mm}

  \ Let us consider also the so-called {\it degenerate} $  \ \Psi \ - \ $ function $ \ \psi_{(r)}(p), \ $ where $  r = \const \in [1,\infty): $

$$
\psi_{(r)}(p) \stackrel{def}{=} 1,  \ p \in [1,r];
$$
so that the correspondent value $  b = b(r) $  is equal to $  r. $  One can  extrapolate formally this function onto the whole  semi-axis $  R^1_+: $

$$
\psi_{(r)}(p)  := \infty, \ p > r.
$$

 \  The classical Lebesgue-Riesz $ L_r  $ norm for the r.v. $  \eta $ is  quite equal to the GLS norm $  ||\eta|| G\psi_{(r)}: $

$$
|\eta|_r = ||\eta|| G\psi_{(r)}.
$$
 \ Thus, the ordinary Lebesgue-Riesz spaces are particular, more precisely, extremal cases of the Grand-Lebesgue ones. \par

\vspace{4mm}

 \ Further, let  $  \phi = \phi(\lambda), \ |\lambda| < \lambda_0 = \const  \in (0, \infty] $
be numerical twice continuous differentiable positive even convex function such that $ \phi(0) = 0; \ \phi(\cdot) $ is
monotonically increasing in the positive semi-interval $  [0, \ \lambda_0) $ and such that
$ \ \phi(\lambda) \sim C(\phi) \cdot \lambda^2, \ \lambda \to 0. $ The set of all such a functions will be denoted by $   \Phi = \{ \ \phi \  \}. $

\vspace{4mm}

 \ {\bf Definition.}  The random variable  $  \ \zeta \ $ belongs to the space  $  \ B(\phi), \ $  for certain
fixed function $  \phi \in \Phi, $ iff there exists a non-negative constant  $ \ \tau \ $ such that

$$
\forall  \lambda: \ |\lambda| < \lambda_0 \ \Rightarrow {\bf E} \exp(\lambda \zeta) \le \exp(\phi(\lambda \tau)). \eqno(2.4)
$$

 \ The minimal value of the constant  $  \ \tau \ $ which satisfies the inequality (2.4) is said to be the
$  B(\phi) \  $ norm  of the r.v. $  \ \zeta: $

$$
||\zeta||B(\phi) \stackrel{def}{=} \max_{\pm} \sup_{\lambda \in (0, \lambda_0)} \phi^{-1} \{ \ln {\bf E} \exp ( \pm \lambda \ \zeta) \}/|\lambda|, \eqno(2.5)
$$
so that

$$
\forall \lambda: |\lambda| < \lambda_0 \ \Rightarrow  {\bf E} \exp(\lambda \zeta) \le \exp ( \phi(\lambda ||\zeta||B(\phi) ) ). \eqno(2.6)
$$

  \ We suppose in fact that the  r.v. $ \ \zeta \ $  is mean zero and satisfies the well-known Cramer's
condition. In this case the generated function $  \ \phi(\cdot) $   may be introduced naturally:

$$
\phi_{\zeta}(\lambda) := \max_{\pm} \ln {\bf E} \exp(\pm \lambda \ \zeta).
$$

 \ See for example the equalities (1.10), (1.11). \par

 \ This natural function  $ \ \phi_{\zeta}(\lambda) \ $ play a very important role in the theory of Large
Deviations (L.D.) Namely, it is well known that

$$
\lim_{n \to \infty} \left[ {\bf P} (S(n)/n > x ) \right]^{1/n} = - \phi^*_{\xi}(x), \ x > 0.
$$

\ Analogous result for the self-normalized sums was obtained by Qi-Man Shao in an article [18]:

$$
\lim_{n \to \infty} \left[ {\bf P} (S(n)/(V(n) \sqrt{n} ) > x ) \right]^{1/n}   =
\sup_{c > 0} \inf_{t > 0} {\bf E} \exp \left[ t( c \xi - x(\xi^2 + c^2)/2    )   \right].
$$

 \ These spaces are complete Banach functional and rearrangement invariant, as
well as considered before Grand Lebesgue Spaces. They were introduced at first
in the article [11]; the detail investigation of these spaces may be found in the
the monographs [2] and [12], chapters 1,2. \par

 \ It is known that if  $  \ \lambda_0 = \infty   $ and $  0 \ne \zeta \in B(\phi) $ if and only if $  \ {\bf E} \zeta = 0 $ and

$$
\exists K = \const \in (0, \infty) \ \Rightarrow \max \left[  {\bf P} (\zeta \ge u), {\bf P}(\zeta \le - u) \right] \le
\exp \left\{ - \phi^*(u/K)   \right\}, \ u > 0,
$$

and herewith

$$
||\zeta||B(\phi) \le C_1(\phi) K \le C_2(\phi) ||\zeta||B(\phi).
$$

 \ More exactly, if $  0 < ||\zeta||B(\phi) = ||\zeta|| < \infty, $ then

$$
\max \left[ {\bf P}(\zeta \ge u), \ {\bf P}(\zeta \le - u)   \right] \le \exp \left(  - \phi^*(u/ ||\zeta||)  \right). \eqno(2.7)
$$

 \ If the r.v. $  \zeta $  belongs to some $ \ B(\phi) \ $ space, then it belongs also to certain  $ \ G\psi \ $ space with

$$
\psi = \psi_{\phi}(p) = \frac{\phi^{-1}(p)}{p}, \ p \ge 1. \eqno(2.8)
$$

 \ The inverse conclusion in not true. Namely, the mean zero r.v. $ \zeta $  can has finite
all the moments $ \ |\zeta|_p < \infty, \ p \ge 1,  \ $ but may not satisfy the Cramer's condition. \par
 \ A very popular class of these spaces form the subgaussian random variables, i.e. for
which  $  \ \phi(\lambda) = \lambda^2 $ and $  \lambda_0 = \infty. $
 The correspondent   $  \ \psi \ $ function has a form  $ \ \psi(p) = \psi_2(p) = \sqrt{p}. $ \par

More generally, suppose

$$
 \phi(\lambda) = \phi_m(\lambda) = |\lambda|^m/m, \  |\lambda| \ge 1, \ \lambda_0 = \infty, \ m = \const > 0. \eqno(2.9)
$$
 \ The correspondent    $ \psi $  function has a form

$$
\psi(p) = \psi_m(p) = p^{1/m}
$$
 and the correspondent tail estimate is follow:

$$
\max \left[  {\bf P}(\zeta \ge u), \ {\bf P} (\zeta \le - u)  \right] \le \exp \left\{ -  (u/K)^m   \right\}, \ u > 0. \eqno(2.10)
$$

 \ These space are used for obtaining of exponential estimates for sums of independent random variables, see e.g.
[11];  [12],  sections 1.6, 2.1-2.5. Indeed, introduce for any function $ \ \phi(\cdot) \ $ from the set $ \ \Phi \ $
a new function   $   \  \overline{\phi}(\cdot)   \  $  which belongs also at the same set:

$$
 \overline{\phi}(\lambda)  \stackrel{def}{=} \sup_{n = 1,2, \ldots} n \  [ \ \phi \{\lambda/\sqrt{n} \} \ ]. \eqno(2.11)
$$
 \ It is easily to see that

$$
\sup_n {\bf E} \exp(\lambda S(n)/\sqrt{n}) \le \exp [ \overline{\phi}(\lambda)] \eqno(2.12)
$$
with correspondent uniform relative the variable $  \  n  \ $ exponential tail  estimate. \par

 \ For instance, if  for some value $  \ m = \const > 0 $

$$
\max \left[  {\bf P}(\xi \ge u), \ {\bf P} (\xi \le - u)  \right] \le \exp \left\{ -  u^m   \right\}, \ u > 0, \eqno(2.13)
$$
then

$$
\max \left[  {\bf P}(S(n)/\sqrt{n} \ge u), \ {\bf P} (S(n)/\sqrt{n} \le - u)  \right] \le
$$

$$
\exp \left\{ - C(m)  u^{\min(m,2)}   \right\}, \ u > 0, \ C(m) \in (0, \infty), \eqno(2.14)
$$
and the last estimate is essentially non-improvable. \par

 \vspace{4mm}

\section{ Main results. Exponential level.}

 \vspace{4mm}

 \  Denote

$$
 \beta_n(\theta; B)  =  \beta(\theta; n, B) :=  n \phi(\theta/\sqrt{n}, \ B \ \theta/n);
$$

$$
\beta^*(z; n, B) := \sup_{\theta > 0} (\theta z - \beta(\theta; n, B)),
$$

\vspace{4mm}

{\bf  Theorem 3.1. }

$$
Q_n(B) \le \exp  \left( - \beta^*(B; n, B ) \right). \eqno(3.0)
$$

\ {\bf Corollary 3.1.}

$$
Q(B) \le \sup_n \exp  \left( - \beta^*(B; n, B) \right). \eqno(3.0a)
$$

\vspace{4mm}

 \ {\bf Proof.}   \ We have:

$$
Q_n(B) = {\bf P} \left(  \sqrt{n} \cdot \frac{\sum \xi(i)}{\sum \xi^2(i)}  \right) =
{\bf P} \left( \sqrt{n} \sum \xi(i) - B \sum \xi^2(i) > 0  \right) =
$$

$$
{\bf P}  \left(  \sum (\sqrt{n} \xi(i) + B( (\sigma^2 - \xi^2(i)  )    ) > n B \sigma^2 \right) =
$$

$$
 {\bf P}  \left(  n^{-1} \sum _i \eta(i,n) >   B \sigma^2 \right). \eqno(3.1)
$$

 \ Since the random variables $   \eta(i)  = \eta(i,n), \ i = 1.2.\ldots,n) $ are centered and common independent,

$$
{\bf E} \exp \left( \theta \cdot n^{-1} \sum \eta(i,n)    \right) = \left[ {\bf E} \exp (\theta \ \eta /n )   \right]^n =
$$

$$
\exp \left(  n \phi(\theta/\sqrt{n}, \ B \ \theta/n)  \right) = \exp \beta(\theta; n, B). \eqno(3.2)
$$

 \ It remains to use the estimate (2.7). \par

 \vspace{4mm}

\section{ Main results. Power level}

 \vspace{4mm}

 \ Let us introduce  the natural $ \  G\psi \ $ function for the r.v., more exactly for the sequence of r.v. $ \ \gamma(i) = \gamma(i,n): $

$$
\Delta(p) = \Delta(p;n) := |\gamma(i,n)|_p = | \ \xi + n^{-1/2} \ B (\sigma^2 - \xi^2)|_p,   \eqno(4.0)
$$
if of course it is finite for some value $  p, \ p \ge 1.$ The sufficient condition for this conclusion is following:

$$
\xi \in \cup _{s = 2}^{\infty} L_s(\Omega, {\bf P})
$$
and does not dependent on the variable $  n. $ \par
 \ We deduce

$$
Q_n(B) = {\bf P}\left (   \sum \eta(i) > n B \sigma^2 \right) =  {\bf P} \left( n^{-1/2} \sum \gamma(i) > B \sigma^2 \right).
\eqno(4.1)
$$

 \ We need to use the famous Rosenthal's inequality, see [20]:

$$
| n^{-1/2} \sum \gamma(i)   |_p \le K_R \cdot \frac{p}{\ln p} \cdot |\eta|_p, \ p \in \supp \ \Delta(\cdot). \eqno(4.2)
$$
where the exact value of the Rosenthal's constant  may be found in [14], [15]:  $ \ K(R) \approx 0.6379\ldots .  $\par

\vspace{3mm}

 \  Introduce  for such a function $  \ \Delta(p) \ $ its Rosenthal transform $  K_R[\Delta](p) $ as follows

$$
K_R[\Delta](p) =  K_R[\Delta; n](p) \stackrel{def}{=}  K_R \cdot \frac{p}{\ln p} \cdot  \Delta(p;n). \eqno(4.3)
$$

 \vspace{4mm}

 \ We obtain by virtue of Rosenthal's inequality

$$
|| n^{-1/2} \sum \gamma(i)   || GK_R[\Delta] \le 1.  \eqno(4.4)
$$

 \ It follows immediately from the estimate(2.2) the next proposition. \par

\vspace{4mm}

{\bf Theorem 4.1.}

\vspace{3mm}

$$
Q_n(B) \le \exp \left\{ - v^*_{ K(R; \Delta, n)} (\ln B)  \right\}, \ B \ge e, \eqno(4.5)
$$
and therefore

$$
Q(B) \le \sup_n \exp \left\{ - v^*_{ K(R; \Delta, n)} (\ln B)  \right\}, \ B \ge e. \eqno(4.6)
$$

\vspace{4mm}

\section { Lower bounds for introduces tail probabilities.}

\vspace{4mm}

 \ Suppose in this section for simplicity  $ \ \sigma = 1.    \ $ It follows immediately from the classical CLT and LLN that

$$
Q(B) \ge \lim_{n \to \infty}  Q_n(B) \ge \exp \left(  - B^2/2  \right), \ B \ge 1. \eqno(5.1)
$$
a first trivial estimate.  A second one is follows:

$$
Q(B) \ge Q_1(B) = {\bf  P} (\xi \ge B \ \xi^2)  = {\bf P}( 0 \le \xi  \le 1/B), \ B > 1; \eqno(5.2)
$$
so that

$$
Q(B) \ge Q_1(B) = \int_{0+}^{(1/B)+} dF_{\xi}(x),
$$
and in the case when the r.v $  \ \xi \ $ has a density of distribution $  f_{\xi}(x), $

$$
Q(B) \ge Q_1(B) = \int_{0+}^{1/B+} f_{\xi}(x) \ dx. \eqno(5.3)
$$

 \ If in addition the density function  $ \ f_{\xi}(x) \   $ is bounded from below in some positive neighborhood of origin, say

$$
f_{\xi} (x)  \ge c_1  = \const > 0, \ x \in (0, 1),
$$
then

$$
Q(B) \ge Q_1(B)  \ge  c_1/B,  \ B \ge 1. \eqno(5.4)
$$

 \ If the function $ f{\xi}(x)  $ is right continuous  at the point $  x = 0+ $ and $  \ f_{\xi}(0) = c_2 = \const > 0, \ $  then obviously

$$
Q_1(B) \sim  c_2/B, \ B \to + \infty.   \eqno(5.5)
$$
 \ On the other words,

$$
Q_1(B) \asymp  c_2/B, \ B \ge 1. \eqno(5.5a)
$$

 \ This state is very different from the other norming function  statement, see  [16]-[20]; as well as
from the alike estimate for one for classical estimates for i., i.d. r.v., see (2.14). \par

\vspace{4mm}

\section{ Examples.}

\vspace{4mm}

 \ {\bf Example 6.0.} Let the r.v. $  \ \xi(i) \ $ be the Rademacher sequence, i.e. they  are independent and

$$
{\bf P} (\xi(i) = 1) = {\bf P} (\xi(i) = - 1) = 1/2.
$$

 \ Then $  \   \xi^2(i) = 1 \  $ and hence $  V(n) = n. $ We return to the classical Khinchine theorem

$$
|\ln Q(B)|  \asymp B^2/2. \ B \in (1, \infty). \eqno(6.1)
$$

\vspace{4mm}

 \ {\bf Example 6.1.} Let now $  \ \{\xi(i) \} $ be independent standard Gaussian distributed r.v. We deduce after some
computations by virtue of theorem 3.1

$$
Q(B) \le c_3 /B, \ B \ge 1; \eqno(6.2)
$$
and the right-hand side of the last inequality coincides up to multiplicative constant with the low estimate (5.5). \par

\begin{center}

 \vspace{5mm}

 {\bf References.}

\end{center}

 \vspace{5mm}

{\bf 1. Bennet C., Sharpley R. } {\it Interpolation of operators.} Orlando, Academic
Press Inc., (1988). \\

\vspace{4mm}

{\bf 2. Buldygin V.V., Kozachenko Yu.V.} {\it Metric Characterization of Random
Variables and Random Processes. } 1998, Translations of Mathematics Monograph, AMS, v.188.\\

\vspace{4mm}

{\bf 3. Caballero M.E., Fernandes B. and Nualart D. }  {\it Estimation of densities
and applications.}  J. of Theoretical Probability, 1998, 27, 537 - 564. \\

\vspace{4mm}

{\bf 4. Xiequan Fan.} {\it Self - normalized deviations with applications to t - statistics.} \\
arXiv 1611.08436 [math. Pr] 25 Nov. 2016. \\

\vspace{4mm}

{\bf 5. A. Fiorenza.} {\it Duality and reflexivity in grand Lebesgue spaces.} Collect.
Math. 51, (2000), 131 - 148.\\

\vspace{4mm}

{\bf 6. A. Fiorenza and G.E. Karadzhov.} {\it Grand and small Lebesgue spaces and
their analogs.} Consiglio Nationale Delle Ricerche, Instituto per le Applicazioni
del Calcoto Mauro Picone”, Sezione di Napoli, Rapporto tecnico 272/03, (2005). \\

\vspace{4mm}

{\bf 7. Gine E., Goetze F. and Mason D. } (1997). {\it When is the Student t -
statistics asymptotically standard normal?} Ann. Probab., 25, (1997), 1514-1531.\\

\vspace{4mm}

{\bf 8. I.Grama, E.Haeusler. }  {\it Large deviations for martingales. } Stochastic Processes and
Applications, 85, (2000), 279-293.\\

\vspace{4mm}

{\bf 9. T. Iwaniec and C. Sbordone.} {\it On the integrability of the Jacobian under
minimal hypotheses.} Arch. Rat.Mech. Anal., 119, (1992), 129-143.\\

\vspace{4mm}

{\bf 10. B.Y.Jing, H.Y.Liang, W.Zhou.}  {\it Self - normalized moderate deviations
for independent random variables.}  Sci. China Math., 55 (11), 2012, 2297-2315.\\

\vspace{4mm}

{\bf 11. Kozachenko Yu. V., Ostrovsky E.I.}  (1985). {\it The Banach Spaces of
random Variables of subgaussian Type.} Theory of Probab. and Math. Stat. (in
Russian). Kiev, KSU, 32, 43-57.\\

\vspace{4mm}

{\bf 12. Ostrovsky E.I.} (1999). {\it Exponential estimations for Random Fields and
its applications,} (in Russian). Moscow-Obninsk, OINPE.\\

\vspace{4mm}

{\bf 13. Ostrovsky E. and Sirota L.}  {\it Vector rearrangement invariant Banach
spaces of random variables with exponential decreasing tails of distributions.} \\
arXiv:1510.04182v1 [math.PR] 14 Oct 2015 \\

\vspace{4mm}

{\bf 14. Ostrovsky E. and Sirota L. }   {\it Schlomilch and Bell series for Bessel's
functions, with probabilistic applications.} \\
arXiv:0804.0089v1 [math.CV] 1 Apr 2008 \\

\vspace{4mm}

{\bf 15. Ostrovsky E. and Sirota L.}  {\it Sharp moment estimates for polynomial martingales.}\\
arXiv:1410.0739v1 [math.PR] 3 Oct 2014 \\

\vspace{4mm}

 {\bf 16. De La Pena V.H. }  {\it A general class of exponential inequalities for martingales and ratios.} 1999,
Ann. Probab., 36, 1902-1938.\\

\vspace{4mm}

{\bf 17. De La Pena V.H., M.J.Klass, T.L.Lai.} {\it Self-normalized Processes:
exponential inequalities, moment Bounds and iterative logarithm law.} 2004 ,
Ann. Probab., 27, 537-564. \\

\vspace{4mm}

{\bf 18. Q.M.Shao.} {\it Self-normalized large deviations.} Ann. Probab., 25, (1997), 285-328.\\

\vspace{4mm}

{\bf 19. Q.M.Shao. } {\it Self-normalized limit theorems: A survey.} Probability
Surveys, V.10 (2013), 69-93. \\

\vspace{4mm}

{\bf 20.  H. P. Rosenthal.} {\it On the subspaces of $ \ L(p) \ (p > 2) \ $  spanned by sequences of
independent random variables. } Israel J. Math., 8  (1970), pp. 273-303.\\

\vspace{4mm}

{\bf 21. Q.Y.Wang, B.Y.Jing.} {\it An exponential non-uniform Berry-Essen for
self - normalized sums.}  Ann. Probab., 27, (4), (1999), 2068-2088.\\

\end{document}